\documentstyle[11pt]{article}
\def\pl{\partial}

\def\*{\raisebox{.5mm}{*}}

\def\div{{\,\rm div\,}}

\def\pa{\partial}

\def\R{I\!\!R}
\def\A{{\cal A}}
\def\H{{\cal H}}

\def\Ga{\Gamma}

\def\<{\langle}
\def\>{\rangle}

\def\Om{\Omega}

\def\det{\mbox{det\,}}

\def\bma{\left[\begin{array}}
\def\ema{\end{array}\right]}
\def\bda{\left|\begin{array}}
\def\eda{\end{array}\right|}
\def\be{\begin{equation}}
\def\ee{\end{equation}}

\newtheorem{thm}{{}\hskip\parindent Theorem}[section]
\newtheorem{lem}{{}\hskip\parindent Lemma}[section]
\newtheorem{pro}{{}\hskip\parindent Proposition}[section]
\newtheorem{exl}{{}\hskip\parindent Example}[section]

\newtheorem{dfn}{{}\hskip\parindent Definition}[section]
\newtheorem{rem}{{}\hskip\parindent Remark}[section]

\large
\title{ Escape Metrics and Its Applications}
\date{}
\author{Zhen-Hu Ning\thanks{Corresponding author, E-mail address: nzh41034@163.com}, Fengyan Yang and Xiaopeng Zhao}
\begin{document}
\maketitle

\footnote{Zhen-Hu Ning,
 Faculty of Information Technology, Beijing University of Technology, Beijing, 100124, China,
 E-mail address: nzh41034@163.com.

\ \ Fengyan Yang,
Key Laboratory of Systems and Control, Institute of Systems Science, Academy of Mathematics and
Systems Science, Chinese Academy of Sciences, Beijing, 100190, China,
E-mail address: yangfengyan12@mails.ucas.edu.cn.

\ \
Xiaopeng Zhao, School of Mathematics, Southeast University, Nanjing, 210096, China;
School of Science, Jiangnan University,
Wuxi,   214122, China,
E-mail address: zhaoxiaopeng@jiangnan.edu.cn.
}

\begin{quote}
\begin{small}

{\bf Abstract} \,\,  Geodesics escape is widely used to study the scattering of
hyperbolic equations. However, there are few progresses  except in a simply connected complete Riemannian manifold with nonpositive curvature.

\quad \ \ We propose a kind of complete Riemannian metrics in $\R^n$, which is called as escape metrics. We expose the relationship between escape metrics and geodesics escape in $\R^n$.
Under the escape metric $g$, we prove that each geodesic of $(\R^n,g)$  escapes, that is, $\lim_{t\rightarrow +\infty} |\gamma (t)|=+\infty$ for any $x\in \R^n$ and any unit-speed  geodesic $\gamma (t)$ starting at $x$. We also obtain the geodesics escape velocity and give the counterexample that if  escape metrics are not satisfied, then there exists an unit-speed geodesic $\gamma (t)$ such that $\overline{\lim}_{t\rightarrow +\infty} |\gamma (t)|<+\infty$.

\quad \ \ In addition,  we establish Morawetz multipliers in Riemannian geometry to derive dispersive estimates for the wave equation on an exterior domain of $\R^n$ with  an escape metric. More concretely, for radial solutions, the uniform decay rate of the local energy is independent of the parity of the dimension $n$. For general solutions, we prove the space-time estimation  of the energy and uniform decay rate $t^{-1}$ of the local energy.  It is worth pointing out that  different from the assumption of  an Euclidean metric at infinity in the existing studies, escape metrics are more general Riemannian metrics.
\\[3mm]
{\bf Keywords}\,\,\,  Escape metrics, Geodesics escape, Wave equation, Dispersive estimates  
\\[3mm]
{\bf AMS(MOS) Subject Classification}\ \ 53C22, 58J45
\end{small}
\end{quote}
\vskip .5cm
\tableofcontents
\def\theequation{1.\arabic{equation}}
\setcounter{equation}{0}
\section{Introduction and Main results}
\vskip .2cm
\subsection{Notations}
\quad  \ \ Let $O$ be the original point of $\R^n$ ($n\ge 2$) and
 \be r(x)= |x|,\quad  x\in\R^n\ee
be  the standard distance function of $\R^n$.
Moreover, let $\<\cdot,\cdot\>$, $\div$, $\nabla$, $\Delta$ and $I_n=(\delta_{i,j})_{n\times n}$ be the standard inner product of $\R^n$,  the standard divergence operator of $\R^n$, the standard gradient operator of $\R^n$, the standard Laplace operator of $\R^n$ and the unit matrix.

Suppose that $(\R^n,g)$ is a smooth complete Riemannian manifold with
\be g=\sum^n_{i,j=1}g_{ij}(x)dx_idx_j,\quad x\in \R^n,\ee such that
  \be \label{1cpde.3} G(x)\frac{\partial}{\partial r} = \frac{\partial}{\partial r} ,\quad |x|\ge r_c,\ee
 where $r_c$ is a positive constant and
\be  G(x)=(g_{ij}(x))_{n\times n},\quad x\in \R^n.\ee

Denote
\be\<X,Y\>_g=\<G(x)X,Y\>,\quad|X|_g^2=\<X,X\>_g,\quad X,Y\in\R^n_x,x\in\R^n.\ee
Let $(r,\theta)$=$(r,\theta_1,\theta_2,\cdots,\theta_{n-1})$ be the polar coordinates of $x\in \R^n$ in the Euclidean metric. Note that
  \be \left\<\frac{\partial}{\partial r},\frac{\partial}{\partial r}\right\>_g(x)=1, \quad \left\<\frac{\partial}{\partial r},\frac{\partial}{\partial\theta_i}\right\>_g(x)=0\quad for \quad  1 \leq i\leq n-1,|x|\ge r_c.\ee
Then
  \be\label{1cpde.2} g=dr^2+\sum_{i,j=1}^{n-1}\gamma_{ij}(r,\theta)d\theta_id\theta_j,\quad |x|\ge r_c,\ee
which  implies $|x| -r_c$ is the  geodesic distance function of $(\R^n,g)$ from $x$ with $|x|\ge r_c$  to $r_c(x/|x|)$.

Let
\be \label{1cpde.4} \Upsilon(x)=(\gamma_{ij})_{(n-1)\times (n-1)}(x),\quad |x|\ge r_c. \ee

Let $S(r)$ be the sphere in $\R^n$ with a radius $r$. Then
\be\left\<X,\frac{\partial}{\partial r}\right\>_g=\left\<X,\frac{\partial}{\partial r}\right\>=0,\quad  \textmd{for }\ \ X\in S(r)_x, |x|\ge r_c.\ee

Let $D$ be the Levi-Civita connection of the metric $g$ and H   a vector field,
then the covariant differential $DH$ of the vector field H is a tensor field of rank 2 as follow:
\be DH(X,Y)(x)=\<D_YH,X\>_g(x)\quad  X,Y\in\R^n_x, x\in\R^n.\ee

Finally, we set $\div_g$, $\nabla_g$ and $\Delta_g$  as the divergence operator of $(\R^n,g)$,  the gradient operator  of $(\R^n,g)$ and the Laplace$-$Beltrami operator of $(\R^n,g)$, repsectively.
\subsection{Escape metrics}

\quad \ \

In this paper, we introduce escape metrics as follows.
\begin{dfn}
We say  $g$  is an {\bf escape metric} if

 \be   \label{2wg.7.4}\left\< \left(\frac{1}{2}\frac{\partial G(x)}{\partial r}\right)X,X \right\>\ge \alpha(x)|X|^2_g\quad for \quad X\in S(r)_x,\ \ |x|\ge r_c, \ee
  \be \label{2wg.7.7} D^2r^2(X,X)\ge 2\rho_c |X|^2_g\quad for \quad X\in \R^n_x,\ \ |x|< r_c, \ee
where $r_c$ is given by (\ref{1cpde.3}), $\rho_c\leq 1$ is a positive constant, $D^2r^2$ is the Hessian of $r^2$ in the metric $g$
and $\alpha(x)$  is a smooth   function defined on $|x|\ge r_c$ satisfying
\be\alpha(x)> -\frac{1}{r},\quad |x|\ge r_c.\ee

\end{dfn}
\begin{rem} If $r_c=0$ in (\ref{1cpde.3}), then  $r(x)=|x|$ is the  geodesic distance function of $(\R^n,g)$ from $x$ to  the original point $O$, (\ref{2wg.7.7}) always  holds true for sufficiently small $r_c$.\end{rem}
\begin{rem}
Escape metrics can be checked by the following Proposition \ref{2wg.9.2}.
\end{rem}

For an exterior domain of $\R^n$, we introduce exterior escape metrics as follows.
\begin{dfn}
We say  $g$ is an {\bf exterior escape metric} if

 \be  \left\< \left(\frac{1}{2}\frac{\partial G(x)}{\partial r}\right)X,X \right\>\ge \alpha(x)|X|^2_g\quad for \quad X\in S(r)_x,\ \ |x|\ge r_c, \ee
where $r_c$ is given by (\ref{1cpde.3}) and $\alpha(x)$  is a smooth   function defined on $|x|\ge r_c$ satisfying
\be\alpha(x)> -\frac{1}{r},\quad |x|\ge r_c.\ee

\end{dfn}

\begin{rem}
Exterior escape metrics can be checked by the following Proposition \ref{2wg.10}.
\end{rem}
\subsection{Geodesics escape}
\quad \ \  One of the important application of  escape metrics is to solve geodesics escape problem.

 Geodesics escape is a very important research topic and there are already some essential researches on it. Bangert\cite{B1981} proved the existence of escaping geodesic without self-intersections in a complete Riemannian manifold homeomorphic to the plane;
 Fernandez and Melian \cite{F2001} studied the quantity of the escaping geodesics of a 2-dimensional, oriented and noncompact Riemann  surface of constant negative curvature;
Gony\cite{G2008} gave the Hausdorff dimension of the terminal points of  escaping geodesics in a hyperbolic manifold.

The assumption that each geodesic of $(\R^n,g)$  escapes, known as non-trapping assumption, that is,
$ \lim_{t\rightarrow +\infty} |\gamma (t)|=+\infty$ for any $x\in \R^n$ and any unit-speed geodesic $\gamma (t)$ starting at $x$,  is widely used to study the scattering of
hyperbolic equation (see for example \cite{w20,1ww,ww1,Burq 2004,Burq 2010,D2009,w3,H2005,H2016,Hidano 2010,76,w32,Zhang 2015} and the references cited therein).  It is well-known that all geodesics escape for a simply connected complete Riemannian manifold with nonpositive curvature. However, for general metrics, since the geodesic is dependent on the nonlinear ordinary differential equation, it is hard to check the non-trapping assumption. 
 Therefore, some effective criterions are further needed to make the non-trapping assumption checkable for general metrics.

 The following two theorems show how geodesics escape under  escape metrics in $\R^n$.
\begin{thm}\label{cmp.3_1}

 Let $g$ be an escape metric such that

 \be r\alpha(x)+1\ge \rho_0,\quad |x|\ge r_c,\ee
 where $\rho_0\leq \rho_c $ is a positive constant and $\rho_c $ is given by (\ref{2wg.7.7}).
Let
\be c_0=\sup_{|x|\leq  r_c}r|Dr|_g(x).\ee
Then, for any $x\in \R^n$ and any unit-speed geodesic $\gamma (t)$ starting at $x$, there exists $c(x)>0$ such that
  \be\label{wg.30_4} |\gamma (t)|\ge \rho_0t -\max\{|x|,c_0\},\quad \forall t>c(x).\ee
 which implies
  \be \lim_{t\rightarrow +\infty} |\gamma (t)|=+\infty.\ee
\end{thm}

\begin{thm}\label{cmp.3_2}

  Let $g$ be an escape metric and


 \be f(y)=\inf_{r_c\leq |x|\leq y} \left\{\alpha(x)+ \frac{1}{|x|}\right\},\quad y\ge r_c.\ee
Then,  for any $x\in \R^n$ and any unit-speed geodesic $\gamma (t)$ starting at $x$, there exists $c(x)>0$ such that
   \begin{eqnarray}\label{cmp.38}
|\gamma (t)| \ge \int_{c(x)}^{t} \left(1-\frac{2}{1+e^{2\int_0^y f(|x_0|+z)dz}}\right)dy +|\gamma (c(x))|,\quad t>c(x).
\end{eqnarray}
 which implies
  \be \lim_{t\rightarrow +\infty} |\gamma (t)|=+\infty.\ee
 Specially, if
\be  \int_{r_c}^{+\infty} f(z)dz =+\infty. \ee
 Then
  \be \label{cmp.39} \lim_{t\rightarrow +\infty}\frac{|\gamma (t)|}{t}=1.\ee

\end{thm}

For  exterior escape metrics, we have similar conclusions with escape metrics.
\begin{thm}\label{cmp.3_3}

  Let $g$ be an exterior escape metric.
Then,  for any $|x|\ge r_c$ and any unit-speed  geodesic $\gamma (t)$ starting at $x$, the following conclusions hold true:
\begin{itemize}\item
If $|x|=r_c$ and $\<\gamma' (t),\frac{\partial}{\partial r}\>_g\ge 0 $, then
\be \lim_{t\rightarrow +\infty} |\gamma (t)|=+\infty.\ee
\item
If $|x|>r_c$, then
\be \lim_{t\rightarrow +\infty} |\gamma (t)|=+\infty,\ee
or  there exist  positive constants $c(x), t_0=t_0(x,\gamma' (0))$ such that
\be 0< t_0\leq c(x)\quad and \quad |\gamma (t_0)|=r_c.\ee
\end{itemize}
\end{thm}
\begin{rem}
If the geodesic is reflected at $|x|=r_c$, then for  exterior escape metrics, for any $|x|\ge r_c$ and any unit-speed geodesic $\gamma (t)$ starting at $x$,
 \be \lim_{t\rightarrow +\infty} |\gamma (t)|=+\infty.\ee
\end{rem}

The following theorem shows the necessity of  escape metrics for geodesics escape.
\begin{thm}\label{cmp.3_4}
Assume that
    \be \frac{\partial G(x)}{\partial r}=- \frac{2}{R_0}\left(I_n-\frac{x\otimes x}{|x|^2}\right)G(x), \quad |x|=R_0,\ee
    where $R_0>r_c$ is a positive constant.
Then, for any $x\in S(R_0)$ and any unit-speed geodesic $\gamma (t)$ starting at $x$ with
 \be\gamma' (0)\in S(R_0)_x,\ee
 we have
  \be\gamma(t)\in S(R_0),\quad \forall t\ge 0.\ee

\end{thm}
\subsection{Wave equation on an exterior domain}

\quad \ \ In the following, we apply exterior escape metrics  to study the dispersive estimates for the wave equation on an exterior domain.

 Let $\Omega$ be an exterior domain in $\R^n$ with a compact smooth boundary
$\Ga$. Assume  that
\be |x|\ge r_c \quad for\ \ any \ \  x\in \Om, \ee
where $r_c$ is given by (\ref{1cpde.3}).
Denote
\be \label{cpde.3.1}r_0= \inf_{x\in \Ga}|x|,\quad r_1= \sup_{x\in \Ga}|x|.\ee
Then $r_1\ge r_0\ge r_c$.

Consider the following system.  
\begin{equation}
\label{wg.1} \cases{u_{tt}-\Delta_g u=0\qquad ~~~~~~~~ (x,t)\in \Omega\times
(0,+\infty),\cr
 u\large |_{\Ga}=0\qquad~~~~~~~~~~~~~~ t\in(0,+\infty),
\cr u(x,0)=u_0(x),\quad u_t(x,0)=u_1(x)\qquad x\in \Omega.}
\end{equation}
Define the energy of the system (\ref{wg.1}) by
\begin{equation}
\label{wg.3}  E(t)=\frac12\int_{\Omega}\left(u_t^2+|\nabla_g u|_g^2\right)dx_g,
\end{equation}
where
\be dx_g=\sqrt{\det(G(x))}dx.\ee
For $a>r_1$, the local energy for the system (\ref{wg.1}) is defined by
\begin{equation}
\label{wg.2}  E(t,a)=\frac12\int_{\Omega(a)}\left(u_t^2+|\nabla_g u|_g^2\right)dx_g,
\end{equation}
where $\Om(a)=\{x|x\in \Om,|x|\leq a\}$.  In this paper, we are interested in how $E(t,a)$ decays.

 If $ g\equiv\sum^n_{i=1}dx^2_i$, the system (\ref{wg.1}) is known as constant coefficient. 
In the case of constant coefficient,  this problem has a long history and a wealth of results were obtained, see for example  
\cite{w25,36,w11,w4,w16,w6,w7,3,6,w5,5} and so on. 

For the general metric $g$, there exist few studies of the system (\ref{wg.1}). To the best of our knowledge, Vodev  \cite{w32} proved uniform local energy decay estimates
on unbounded Riemannian manifold with nontrapping metrics assumption. In \cite{w31}, Christianson obtained a sub-exponential local energy decay estimate on a complete, non-compact, odd-dimensional Riemannian manifold with one trapped hyperbolic geodesic and a Euclidean  metric outside a compact set. Later, Christianson \cite{C2009} derived a polynomial local energy  decay rate on a non-compact Riemannian manifold with a hyperbolic thin trapped set. Vasy\cite{V2013} considered the weighted energy decay for the wave equation at low frequency on scattering manifolds. Besides, we refer the reader to \cite{Blair,ww1,Metcalfe,76,Zhang 2015} for a series of Strichartz estimates for wave equations.
Due to the limit that the metric $g$ is assumed as an Euclidean metric at infinity, we are motivated to study the dispersive estimates for the system (\ref{wg.1})  under a general Riemannian metric.

As is known, the multiplier method is a simple and effective tool for energy  estimate on PDEs. In particular,  the celebrated Morawetz's multipliers introduced by \cite{w6} have been extendedly used for studying the energy decay of the wave equation with
constant coefficients, see \cite{w25,w11,w16,w7,5} and many others. Therefore,
we establish  Morawetz multipliers in Riemannian geometry to derive dispersive estimates for the system (\ref{wg.1}).

Define
\be \widehat{H}^1_{0}(\Om)=\left\{  w(x)  \Big| \ \  w\large|_\Ga=0, \quad \int_{\Omega}|\nabla_g w|_g^2dx_g<+\infty \right\}.\label{wg.4.1}\ee
It is well-known that the  system (\ref{wg.1}) is well-posed with
\be  u_t\in C([0,+\infty),L^2(\Om)) \ \ and \ \ u\in C([0,+\infty),\widehat{H}^1_{0}(\Om)).\ee
Then
\be E(t)=E(0),\quad \forall t>0.\ee

The finite speed of propagation property of the wave  can be stated as follows.

\begin{thm}\label{cmp.3_5} Let initial datum $(u_0,u_1)$ satisfy
\be   u_0(x)=u_1(x)=0, \quad |x|\ge R_0.\ee
where $R_0>r_1$ is a  constant. Then
\be u(x,t)=0\quad \textmd{for}~|x|\ge  R_0+t.\ee
\end{thm}

The following Assumption\ {\bf(A)} and Assumption\ {\bf(B)} are the main assumptions for the wave equation.

{\bf Assumption (A)}\,\,\,  Let $g$ be an  exterior escape metric such that
\be  \label{cmp.1}  \alpha(x)= \frac{ m_1 -1}{r}  \quad for \ \ |x|\ge r_0,\ee
\be   \label{cmp.2} \det(G(x))=\eta_0(\theta)r^{2m_2-2(n-1)}\quad for \ \ |x|\ge r_0.\ee


where   $m_1,m_2$ are positive constants, $\eta_0(\theta)$ is a positive function.
\begin{rem}
 From the following formulas (\ref{cmp.41}) and  (\ref{cmp.61}), we have
      \begin{eqnarray} \frac{m_2}{ r}&&=\frac{n-1}{r}+\frac{\partial\ln \sqrt{\det(G(x))}}{\partial r}=\Delta_g r =tr D^2r\nonumber\\
&&\ge (n-1)\left(\alpha(x)+\frac{1}{r}\right)= \frac{(n-1) m_1 }{r},\quad |x|\ge r_0.     \end{eqnarray}
Then
   \be   m_2\ge (n-1)m_1.\ee
\end{rem}
\begin{exl}
Let $r_0= r_c$ and  $g$  satisfy
\be G(x)=\frac{x\otimes x}{|x|^2}+r^{2(m_1-1)}\left(I_n-\frac{x\otimes x}{|x|^2}\right), \quad   |x|\ge r_0.\ee
Then
  \be  G(x)\frac{\partial}{\partial r} = \frac{\partial}{\partial r} ,\quad |x|\ge r_0,\ee
 \be  \left\< \left(\frac{1}{2}\frac{\partial G(x)}{\partial r}\right)X,X \right\>= \frac{ m_1 -1}{r} |X|^2_g\quad for \quad X\in S(r)_x,\ \ |x|\ge r_0, \ee
\be  \det(G(x))=r^{2 (m_1-1)(n-1)} \quad for \quad   |x|\ge r_0.\ee
\end{exl}

{\bf Assumption (B)}\,\,\, Let $g$ be an  exterior escape metric such that
\be  \label{cmp.6}  \alpha(x)= m_1 r^{-s_1}-\frac{ 1 }{r} \quad for \ \ |x|\ge r_0,\ee
\be   \label{cmp.7} \det(G(x))=\eta_0(\theta)r^{-2(n-1)}e^{2m_2\int_{r_c}^ry^{-s_2}dy} \quad for \ \ |x|\ge r_0,\ee
where $\eta_0(\theta)$ is a positive function, $m_1,m_2,s_1,s_2$ are positive constants such that
\be \label{2wg.9.8}s_1>1,\quad 0<s_2\leq 1 \quad and \quad (s_2+1)r_0^{s_2-1}<m_2.\ee

\begin{rem}
 From the following formulas (\ref{cmp.41}) and  (\ref{cmp.61}), we have
      \begin{eqnarray} m_2 r^{-s_2}&&=\frac{n-1}{r}+\frac{\partial\ln \sqrt{\det(G(x))}}{\partial r}=\Delta_g r =tr D^2r\nonumber\\
&&\ge (n-1)\left(\alpha(x)+\frac{1}{r}\right)=(n-1) m_1 r^{-s_1},\quad |x|\ge r_0.     \end{eqnarray}
Then
   \be   m_2\ge (n-1)m_1r^{(s_2-s_1)}_0.\ee
\end{rem}

\begin{exl}
Let $r_0= r_c$ and  $g$  satisfy
\be  G(x)=\frac{x\otimes x}{|x|^2}+e^{\frac{2\int_{r_0}^rm_2 y^{-s_2}dy-2(n-1)}{n-1}}\left(I_n -\frac{x\otimes x}{|x|^2}\right),\ \ |x|\ge r_0,\ee
where $s_1,s_2,m_2$ are given by (\ref{2wg.9.8}). Then
  \be  G(x)\frac{\partial}{\partial r} = \frac{\partial}{\partial r} ,\quad |x|\ge r_0,\ee
\begin{eqnarray}  \left\< \left(\frac{1}{2}\frac{\partial G(x)}{\partial r}\right)X,X \right\>&&=\left(\frac{m_2 r^{-s_2}}{n-1}-1\right)|X|^2_g \nonumber\\
&& \ge  \left(m_1 r^{-s_1}-\frac{ 1 }{r} \right)|X|^2_g\quad for \quad X\in S(r)_x,\ \ |x|\ge r_0,\end{eqnarray}
 \be   \det(G(x))=r^{-2(n-1)}e^{2m_2\int_{r_c}^ry^{-s_2}dy} \quad for \ \ |x|\ge r_0,\ee
where $m_1$ is a positive constant such that
   \be   m_2\ge (n-1)m_1r^{(s_2-s_1)}_0.\ee

\end{exl}

\begin{rem}
More examples of Assumption\ {\bf(A)} and Assumption\ {\bf(B)} can be obtained by solving $det(G(x))$ from the following (\ref{2wg.9.7}) with a given $\alpha(x)$.
\end{rem}

  In order to facilitate the discussion, we define
   \be \nu(x) \ \ is\ \  the\ \   unit\ \  normal \ \ vector\ \ outside\ \  \Om  \ \ in \ \ (\R^n,g)\ \  for\ \  x\in  \Ga, \ee and
    \be \nu(x)=\frac{\partial}{\partial r}\ \  for\ \  |x|> r_1. \ee

\begin{thm}\label{cmp.3_6} Suppose that the following three conditions hold:
\be \label{wg.7_1}Assumption\ {\bf(A)}\ \ \ holds\ \ true\ \ with \ \ \ m_1>\frac{1}{2}, \ee
\be \label{wg.7_2}\frac{\partial r}{\partial \nu}\le 0,\quad x\in \Ga,\ee
\be \label{wg.7_4}u_0(x)=u_1(x)=0,\quad |x|\ge R_0, \ee
where $R_0>r_1$ is a  constant.

Then
there
exists positive constant $C(a,R_0)$ such that
\begin{equation}
\label{wg.7.1}
E(a,t)\leq \frac{C(a,R_0)}{t} E(0),\quad \forall t> 0.\end{equation}
\end{thm}

\begin{thm}\label{cmp.3_7}  Suppose that the following three conditions hold:
\be \label{cmp.9}Assumption\ {\bf(B)}\ \ \ holds\ \ true, \ee
\be \label{cmp.10}\frac{\partial r}{\partial \nu}\le 0,\quad x\in \Ga,\ee
\be \label{cmp.11}u_0(x)=u_1(x)=0,\quad |x|\ge R_0, \ee
where $R_0>r_1$ is a  constant.

Then
there
exists positive constant $C(R_0)$ such that
\begin{equation}
\label{cmp.12}
\int_0^{+\infty}\int_{\Omega}r^{-s_1}\left(u_t^2+|\nabla_g u|_g^2\right)dx_gdt\leq C(R_0) E(0).\end{equation}
\end{thm}

The organization of the rest of our paper goes as follows. In Section 2, we will  show how to check
 escape metrics and exterior escape metrics. Then   some multiplier identities and key lammas for our problem  will be given in Section 3.  The technical details of the proofs of the results for geodesics escape will be provided in Section 4. Finally, Section 5 is devoted to the wave system (\ref{wg.1}), in which the propagation property, the uniform local energy decay for radial solutions, the uniform local energy decay and the space-time energy estimation for general solutions are present in sequence.
\vskip .5cm
\def\theequation{2.\arabic{equation}}
\setcounter{equation}{0}
\section{Checking the Metric}
\vskip .2cm

\quad Escape metrics can be checked by the following proposition.
\begin{pro}  \label{2wg.9.2} Let
   \begin{eqnarray} G(x)&&=W(x)+e^{\int_{r_c}^r2\alpha(y,\theta)dy}(I_n-W(r_c,\theta))\nonumber\\
&&\quad +e^{\int_{r_c}^r2\alpha(y,\theta)dy}\int_{r_c}^r2e^{-\int_{r_c}^y2\alpha(z,\theta)dz} Q(y,\theta)(I_n-W(y,\theta))dy, \quad   x\in \R^n ,\end{eqnarray}
where
 \be W(x)=\frac{x\otimes x}{|x|^2}, x\in \R^n\backslash O,\ee
$\alpha(x)$ is a smooth   function defined on $\R^n$ satisfying
\be\alpha(x)> -\frac{1}{r},\quad  |x|\ge r_c,\ee
     \be\alpha(x)=0,\quad  |x|< r_c,\ee
and  $Q(x)=(q_{i,j})_{(n-1)\times (n-1)}(x)$ is a smooth, symmetric, nonnegative definite   matrix  function  defined on $\R^n$  satisfying
     \be Q(x)=0,\quad |x|< r_c .\ee

     Then

\be \label{2wg.9.3} G(x)\frac{\partial}{\partial r}=\frac{\partial}{\partial r},\quad x\in\R^n,\ee

    \be   \label{2wg.9.4}D^2r^2(X,X)=2 |X|^2_g\quad for \quad X\in \R^n_x,\ \ |x|< r_c, \ee

     \be  \label{2wg.9.5} \left\< \left(\frac{1}{2}\frac{\partial G(x)}{\partial r}\right)X,X \right\>\ge \alpha(x)|X|^2_g\quad for \quad X\in S(r)_x,\ \ |x|\ge r_c. \ee
      \end{pro}

     {\bf Proof.} Note that
     \be \label{2wg.9.6}  G(x)=I_n ,\quad  |x|<r_c.\ee
     Then   \be D^2r^2(X,X)= 2|X|^2_g\quad for \quad X\in \R^n_x,|x|<  r_c.\ee
     the equality  (\ref{2wg.9.4}) holds.

     Note that
     \be \frac{\partial}{\partial r} =\frac{x}{|x|}.\ee
 Then
         \be W(x)\frac{\partial}{\partial r}=\frac{\partial}{\partial r},\quad  x\in  \R^n\backslash O.\ee
Hence
         \be G(x)\frac{\partial}{\partial r}=\frac{\partial}{\partial r},\quad  x\in  \R^n\backslash O.\ee
With (\ref{2wg.9.6}),  the equality  (\ref{2wg.9.3}) holds.

 Note that  \be W(x)X=0\quad for \quad X\in S(r)_x,x\in R^n\backslash O,\ee and\be  \frac{1}{2}\frac{\partial G(x)}{\partial r} =\left(\alpha(x)G(x)+Q(x)\right)(I_n-W(x)),\quad |x|\ge r_c.\ee

      Then  \begin{eqnarray}  \left\< \left(\frac{1}{2}\frac{\partial G(x)}{\partial r}\right)X,X \right\>&&= \alpha(x)|X|^2_g+\<Q(x)X,X\>\nonumber\\&&\ge \alpha(x)|X|^2_g\quad for \quad X\in S(r)_x,\ \ |x|\ge r_c. \end{eqnarray}
       The inequality (\ref{2wg.9.5}) holds true.$\Box$

By a similar proof with Proposition \ref{2wg.9.2},  exterior escape metrics can checked by the following proposition.
\begin{pro}  \label{2wg.10}  Let
   \begin{eqnarray}\label{2wg.9.7} G(x)&&=W(x)+e^{\int_{r_c}^r2\alpha(y,\theta)dy}P(r_c,\theta)(I_n-W(r_c,\theta))\nonumber\\
&&\quad +e^{\int_{r_c}^r2\alpha(y,\theta)dy}\int_{r_c}^r2e^{-\int_{r_c}^y2\alpha(z,\theta)dz} Q(y,\theta)(I_n-W(y,\theta))dy, \quad   |x|\ge r_c,\end{eqnarray}
where
 \be W(x)=\frac{x\otimes x}{|x|^2}, x\in \R^n\backslash O,\ee
$\alpha(x)$ is a smooth   function defined on $|x|\ge r_c$ satisfying
\be\alpha(x)> -\frac{1}{r},\quad  |x|\ge r_c,\ee
 $P(x)=(p_{i,j})_{(n-1)\times (n-1)}(x)$ is a smooth, symmetric, positive definite   matrix  function  defined on $|x|=r_c$ and  $Q(x)=(q_{i,j})_{(n-1)\times (n-1)}(x)$ is a smooth, symmetric, nonnegative definite   matrix  function  defined on $|x|\ge r_c$.

     Then
     \be G(x)\frac{\partial}{\partial r}=\frac{\partial}{\partial r},\quad |x|\ge r_c,\ee
      \be  \left\< \left(\frac{1}{2}\frac{\partial G(x)}{\partial r}\right)X,X \right\>\ge \alpha(x)|X|^2_g\quad for \quad X\in S(r)_x,\ \ |x|\ge r_c.\ee
      \end{pro}

 \vskip .5cm
\def\theequation{3.\arabic{equation}}
\setcounter{equation}{0}
\section{Multiplier Identities and Key Lammas}
\vskip .2cm

\begin{lem}\label{wg.14}
Suppose that $u(x,t)$ solves the system (\ref{wg.1})
and  $\H$ is a  $C^1$ vector
field defined on $\overline \Om$. Then
\begin{eqnarray}
 \label{wg.14.1}
 &&\int_0^T\int_{\partial\Om(a)}\frac{\pa u}{\pa\nu}\H(u) d\Ga_g dt+\frac12\int_0^T\int_{\partial\Om(a)}
\left(u_t^2-\left|\nabla_g u\right|_g^2\right)\left<\H,\nu\right>_g d\Ga_g dt\nonumber\\
=&&(u_t,\H(u))\Big |^T_0+\int_0^T\int_{\Om(a)}D\H(\nabla_g
u,\nabla_g u) dx_g dt\nonumber\\
&&+\frac12\int_0^T\int_{\Om(a)}\left(u_t^2-\left|\nabla_g u\right|_g^2\right)\div_g\H dx_g dt,
\end{eqnarray}
where
\be (u_t,\H(u))\Big |^T_0 =\int_{\Omega(a)}u_t\H(u) dx_g\Big |^T_0, ~~d\Ga_g=\sqrt{\det(G(x))}d\Ga.\ee

Moreover, assume that $P\in C^2(\overline{\Om})$ and $Q\in C^1(\overline{\Om}\times [0,+\infty ))$ . Then
\begin{eqnarray}
\label{wg.14.2}
\int_0^T\int_{\Om(a)}\left(u_t^2-\left|\nabla_g u\right|_g^2\right)P dx_g dt =&&\frac12\int_0^T\int_{\partial\Om(a)}u^2\frac{\pa
P}{\pa\nu}d\Ga_g dt \int_0^T
-\frac12\int_{\Om(a)}u^2\Delta_gP dx_g dt\nonumber\\
&&-\int_0^T\int_{\partial\Om(a)}Pu\frac{\pa
u}{\pa\nu}d\Ga_g dt+(u_t,u P)\Big|^T_0.\end{eqnarray}

and

\begin{eqnarray}
\label{wg.14.3}
 \int_{\Om(a)}\left(u_t^2+\left |\nabla_g
u\right|_g^2\right)Q dx_g \Big|^T_0=&&\int_0^T
\int_{\Om(a)}\left(u_t^2+\left |\nabla_g
u\right|_g^2\right)Q_t dx_gdt\nonumber\\&&-2\int_0^T
\int_{\Om(a)}u_t\nabla_gQ(u) dx_g dt\nonumber\\
&&+2\int_0^T\int_{\partial\Om(a)}Qu_t\frac{\pa
u}{\pa\nu }d\Ga_g dt.\end{eqnarray}

\end{lem}

{\bf Proof}.
Firstly, we multiply the wave equation in (\ref{wg.1}) by $\H(u)$ and integrate over $\Omega(a)\times
(0,T) $, noting that
\begin{eqnarray} \<\nabla_g f,\nabla_g (H(f))\>_g =&&\nabla_g f\<\nabla_g f,H\>_g
= D^2f(H,\nabla_g f)+DH(\nabla_g f,\nabla_g f)\nonumber\\
=&& D^2f(\nabla_g f,H)+DH(\nabla_g f,\nabla_g f)\nonumber\\
=&&\frac{1}{2}H(|\nabla_g f|_g^2)+DH(\nabla_g f,\nabla_g f)\nonumber\\
=&&DH(\nabla_g f,\nabla_g f)+\frac{1}{2}\div_g(|\nabla_g f|_g^2H)-\frac{1}{2}|\nabla_g f|_g^2\div_g H.\end{eqnarray}
the equality  (\ref{wg.14.1}) follows from Green's formula.

Secondly, we multiply the wave equation in (\ref{wg.1}) by $Pu$ and integrate over $\Omega(a)\times
(0,T)$. The equality  (\ref{wg.14.2}) follows from Green's formula.
Finally, the equality  (\ref{wg.14.3}) follows from Green's formula.$\Box$

 \begin{lem}\label{pde.4}
\be \label{pde.5} D^2r(X,X)=\left\<\left(\frac{1}{2}\frac{\partial G(x)}{\partial r}\right)X,X\right\> +\frac{1}{r}|X|^2_g\quad for\ \  X\in S(r)_x, |x|\ge r_c, \ee
 where $D^2r$ is the Hessian of $r$ in the metric $g$.
  \end{lem}

  {\bf Proof}. Let $\Upsilon(x)$ be given by (\ref{1cpde.4}). Denote
  \be \aleph=\left\{x\ \ \Big|\quad |x|\ge r_c, \det(\Upsilon(x))\neq 0\right\}.\ee
 Let $x\in \aleph$ and denote $\theta_n=r$, we have
   \be \gamma_{ni}(x)=\gamma_{in}(x)=\left\<\frac{\partial}{\partial\theta_i},\frac{\partial}{\partial\theta_n}\right\>_g(x)=0\ \ 1\le i\le n-1,\quad \gamma_{nn}(x)=1.\ee
 Let  $(\gamma^{ij})_{n\times n}(x)=(\gamma_{ij})^{-1}_{n\times n}(x)$.  Then \be \gamma^{ni}(x)=\gamma^{in}(x)=0\ \ 1\le i\le n-1,\quad \gamma^{nn}(x)=1.\ee

 We compute Christofell symbols as
  \be \Gamma^k_{in}=\frac{1}{2}\sum_{l=1}^n\gamma^{kl}\left(\frac{\partial(\gamma_{il})}{\partial r}+\frac{\partial(\gamma_{nl})}{\partial\theta_i}-\frac{\partial(\gamma_{in})}{\partial\theta_l}\right)
  =\frac{1}{2}\sum_{l=1}^{n-1}\gamma^{kl}\frac{\partial(\gamma_{il})}{\partial r},\ee
  for $ 1\le i,k\le n-1$,  which give
  \be D_{\frac{\partial}{\partial \theta_i}}\frac{\partial}{\partial r}=\frac{1}{2}\sum_{k=1}^{n-1}\left(\sum_{l=1}^{n-1}\gamma^{kl}\frac{\partial(\gamma_{il})}{\partial r}\right)\frac{\partial}{\partial \theta_k}.\ee
Then for $X=\sum_{i=1}^{n-1}X_i\frac{\partial}{\partial\theta_i} \in S(r)_x$, we deduce that
\begin{eqnarray} \label{cmp.50} D^2r(X,X)&&= \sum_{i,j=1}^{n-1}\left\< D_{\frac{\partial}{\partial \theta_i}}\frac{\partial}{\partial r},\frac{\partial}{\partial \theta_j}\right\>_gX_iX_j
 \nonumber\\
&&  =\frac{1}{2}\sum_{i,j,k,l=1}^{n-1}\gamma^{kl}\frac{\partial(\gamma_{il})}{\partial r}\gamma_{kj}X_iX_j \nonumber\\
&& =\frac{1}{2}\sum_{i,j=1}^{n-1}\frac{\partial(\gamma_{ij})}{\partial r}X_iX_j  .\end{eqnarray}

Note that
\be \gamma_{ij}(x)=\left\<G(x)\frac{\partial}{\partial \theta_i},\frac{\partial}{\partial \theta_j}\right\>,\qquad \nabla_{\frac{\partial}{\partial r}}\frac{\partial}{\partial \theta_i}(x)=\frac{1}{r}\frac{\partial}{\partial \theta_i}(x), 1\leq i\leq n-1. \ee
Then
\begin{eqnarray}
\frac{\partial \gamma_{ij}}{\partial r}&& =\left\<\frac{\partial G(x)}{\partial r}\frac{\partial}{\partial \theta_i},\frac{\partial}{\partial \theta_j}\right\>+\left\< G(x)\left(\nabla_{\frac{\partial}{\partial r}}\frac{\partial}{\partial \theta_i}\right),\frac{\partial}{\partial \theta_j}\right\>+\left\< G(x)\frac{\partial}{\partial \theta_i},\nabla_{\frac{\partial}{\partial r}}\frac{\partial}{\partial \theta_j}\right\>\nonumber\\
&& =\left\<\frac{\partial G(x)}{\partial r}\frac{\partial}{\partial \theta_i},\frac{\partial}{\partial \theta_j}\right\>+\frac{2}{r}\left\< G(x)\frac{\partial}{\partial \theta_i},\frac{\partial}{\partial \theta_j}\right\>.
\end{eqnarray}
It follows from (\ref{cmp.50}) that
\begin{eqnarray}  D^2r(X,X)&& =\frac{1}{2 }\sum_{i,j=1}^{n-1}\frac{\partial(\gamma_{ij})}{\partial r}X_iX_j \nonumber\\
&& =\left\<\left(\frac{1}{2}\frac{\partial G(x)}{\partial r}\right)X,X\right\> +\frac{1}{r}|X|^2_g.\end{eqnarray}

Note that
\be \aleph \ \  is\ \  dense\ \  in \ \ |x|\ge r_c. \ee
The equality (\ref{pde.5}) holds true. $\Box$

      \begin{lem} \label{cmp.40}
    Let $g$ be an escape metric or an exterior escape metric.
Then
 \be \label{cmp.41} D^2 r(X,X)\ge   \left(\alpha(x)+\frac{1}{r}\right)|X|^2_g\quad \textmd{for }\ \ X\in S(r)_x,|x|\ge r_c.\ee
\end{lem}

  {\bf Proof.}  \ \ \

    With (\ref{pde.5}), we obtain
    \begin{eqnarray}   D^2r(X,X)&&\ge  \left(\alpha(x)+\frac{1}{r}\right)  |X|^2_g \quad \textmd{for }\ \ X\in S(r)_x,|x|\ge r_c.\end{eqnarray}
The inequality  (\ref{cmp.41}) holds true.$\Box$

          \begin{lem} \label{cmp.4}
   Let $g$ be escape metric.
Then
 \be \label{2wg.8} D^2r^2(X,X)\ge  h(x)|X|^2_g\quad \textmd{for }\ \  X\in \R^n_x,x\in \R^n,\ee
 where  \be h(x)= \cases{2\rho_c,\quad |x|< r_c,\cr
     \min\{2,2r\alpha(x)+2\},\quad |x|\ge r_c.}\ee
\end{lem}

  {\bf Proof.}  \ \ \
  With (\ref{pde.5}), we obtain
       \begin{eqnarray}   D^2r(X,X)&&\ge  \left(\alpha(r)+\frac{1}{r}\right)  |X|^2_g \quad \textmd{for }\ \ X\in S(r)_x,|x|\ge r_c.\end{eqnarray}
    Note that
    \be D^2r^2=2Dr\otimes Dr +2r D^2r.\ee
  Then
     \be  D^2r^2(X,X)\ge \min\{2,2r\alpha(x)+2\}|X|^2_g \quad \textmd{for }\ \ X\in \R^n_x, |x|\ge r_c.\ee


With (\ref{2wg.7.7}), the inequality  (\ref{2wg.8}) holds true.$\Box$

    \begin{lem} \label{cmp.60}
    Let $g$ be an escape metric or an exterior escape metric.
Then
 \be \label{cmp.61} \Delta_g r=\frac{n-1}{r}+\frac{\partial\ln \sqrt{\det(G(x))}}{\partial r}, \quad |x|\ge r_c.\ee
\end{lem}

  {\bf Proof.}  \ \ \
With  (\ref{1cpde.3}), for $|x|\ge r_c$, we have
     \begin{eqnarray}   \Delta_g r&& =\div_g\nabla_g r=\frac{1}{\sqrt{\det(G(x))}}\div (\sqrt{\det(G(x))} \frac{\partial}{\partial r})
     \nonumber\\&&= \frac{n-1}{r}+\frac{\partial\ln \sqrt{\det(G(x))}}{\partial r}.\end{eqnarray}

The equality  (\ref{cmp.61}) holds true.$\Box$

\vskip .5cm
\def\theequation{4.\arabic{equation}}
\setcounter{equation}{0}
\section{Proofs for  Geodesics Escape}
\vskip .2cm

\begin{lem}\label{cmp.35}  Let $\widehat{\Omega}$ be a bounded domain of $\R^n$. Assume that there exists a $C^1$ vector field $H$ on $\R^n$ satisfying \be D H(X,X)\geq
\rho_0|X|_g^2\quad for\ \ all\ \ X\in\R^n_x,\quad
x\in\overline{\widehat{\Omega}},\ee
where $\rho_0$ is a positive constant.
Then,  for any $x\in \widehat{\Omega}$ and any unit-speed geodesic $\gamma (t)$ starting at $x$,
if
   \be \gamma (t)\in \widehat{\Omega},\quad 0\leq t \leq t_0, \ee
then
\be\label{cmp.36} t_0 \leq \frac{\sup\left\{|H|_g(x)\Big|\ \  x\in \overline{\widehat{\Omega}}\right\}+|H|_g(x)}{\rho_0}.\ee
\end{lem}

{\bf Proof.} Note that
\be |\gamma' (t)|_g=1,\ \  D_{\gamma' (t)}\gamma' (t)=0 \qquad \forall t\ge 0. \ee Then
\be\<H,\gamma' (t)\>_g\Big|_0^{t_0}= \int_0^{t_0}\gamma' (t)\<H,\gamma' (t)\>_gdt=\int_0^{t_0}D H(\gamma' (t),\gamma' (t))dt\ge \rho_0 t_0. \ee
We have
  \be t_0 \leq\frac{\sup\left\{|H|_g(x)\Big|\ \  x\in \overline{\widehat{\Omega}}\right\}+|H|_g(x)}{\rho_0}.\ee $\Box$

\begin{lem} \label{cmp.5}
Assume that \be \label{wg.30_1}D^2 r^2(X,X)\geq
2\rho_0|X|_g^2\quad for\ \ all\ \ \quad X\in\R^n_x, x\in \R^n,
\ee
where $\rho_0\leq 1$ is a positive constant.
Let
\be c_0=\sup_{|x|\leq  r_c}r|Dr|_g(x).\ee
Then, for any $x\in \R^n$ and any unit-speed geodesic $\gamma (t)$ starting at $x$, there exists $c(x)>0$ such that
  \be\label{wg.30_5} |\gamma (t)|\ge \rho_0t -\max\{|x|,c_0\},\quad \forall t>c(x).\ee
  which implies
  \be \lim_{t\rightarrow +\infty} |\gamma (t)|=+\infty.\ee
\end{lem}

{\bf Proof.}
Note that
\be r|Dr|_g(x)=r(x),\quad |x| \ge r_c.\ee
Hence
\be r|Dr|_g(x)\leq \max\{r(x),c_0\},\quad x\in \R^n.\ee

Let $H=Dr^2$, it is easy to see that
 \begin{eqnarray}
 \<H,\gamma' (z)\>_g\Big|_0^{t}&&= \int_0^{t}\gamma' (z)\<H,\gamma' (z)\>_gdz=\int_0^{t}D H(\gamma' (z),\gamma' (z))dz\ge 2\rho_0t.
\end{eqnarray}
Then
 \be  2\rho_0t \leq 2 \max\{|x|,c_0\} +2\max\{|\gamma (t)|,c_0\}.\ee
Let \be c(x)= \frac{\max\{|x|,c_0\} +c_0}{\rho_0},\ee
 the estimate (\ref{wg.30_5}) holds .
$\Box$

\begin{lem}\label{cmp.31}

 Let $\rho_0\leq 1$ be a positive constant. Assume that
    \be \label{cmp.32}D^2r^2(X,X) \ge 2\rho_0|X|^2_g,\quad X\in\R^n_x,|x|< r_c,\ee
 \be \label{cmp.33}D^2 r(X,X)\geq
f(r)|X|_g^2\quad  \textmd{for }\ \ X\in S(r)_x,|x|\ge r_c,\ee
where  $f\in C([r_c,+\infty))$ is a decreasing positive function.

Then,  for any $x\in \R^n$ and any unit-speed geodesic $\gamma (t)$ starting at $x$, there exists $c(x)>0$ such that
   \begin{eqnarray}\label{cmp.34}
|\gamma (t)| \ge \int_{c(x)}^{t} \left(1-\frac{2}{1+e^{2\int_0^y f(|x_0|+z)dz}}\right)dy +|\gamma (c(x))|,\quad t>c(x),
\end{eqnarray}
 which implies
  \be \lim_{t\rightarrow +\infty} |\gamma (t)|=+\infty.\ee

 Specially, if
\be  \int_{r_c}^{+\infty} f(z)dz =+\infty. \ee
 Then
  \be \label{cmp.30} \lim_{t\rightarrow +\infty}\frac{|\gamma (t)|}{t}=1.\ee

\end{lem}

{\bf Proof.} 
Let $x$ be a fixed point.  From Lemma \ref{cmp.4} and Lemma \ref{cmp.35}, there exists $c(x)>0$, for any unit-speed geodesic $\gamma (t)$ starting at $x$, there exist $0 \leq t_1,t_2 \leq  c(x)$ such that
\be |\gamma (t_1)|=\max\{|x|,r_c\}+\frac{1}{2},\quad |\gamma (t_2)|=\max\{|x|,r_c\}+\frac{3}{2}.\ee

Let $t_0$ satisfy
\be t_0= \sup \left\{\ \ t\ \ \Big| \quad t_1\leq t\leq t_2,\ \    |\gamma (t)|=\max\{|x|,r_c\}+1 \right\}.\ee
Then
  \be |\gamma(t_0)|_t =\left\<\frac{\partial}{\partial r},\gamma' (t_0)\right\>_g \ge 0.\ee

Let \be h(t)=\left\<\frac{\partial}{\partial r}, \gamma' (t)\right\>_g,\ee we have $-1\leq h\leq 1$, $h(t_0)\ge 0$.
For $t\ge t_0$, we deduce that
 \begin{eqnarray}\label{cmp.42}
h_t&&= \gamma' (t)\left\<\frac{\partial}{\partial r}, \gamma' (t)\right\>_g= D^2r(\gamma' (t), \gamma' (t))= D^2r\left(\gamma' (t)-h\frac{\partial}{\partial r}, \gamma' (t)-h\frac{\partial}{\partial r}\right) \nonumber\\
&&\ge f(|\gamma (t)|) (1-h^2).
\end{eqnarray}
which implies $h_t\ge 0$, $h$  is increasing on $ [t_0,+\infty)$. Then
\be h(t)\ge h(t_0)\ge0,\quad t> t_0.\ee
\be |\gamma (t)|=|\gamma (t_0)|+\int_{t_0}^th(z)dz\ge|\gamma (t_0)| ,\quad t>t_0.\ee

Note that $f$ is decreasing and
\be|\gamma (t)|\leq  |x_0|+t.\ee
It follows from (\ref{cmp.42}) that    \be h_t \ge f(|x_0|+t) (1-h^2),\quad t\ge t_0.\ee
Then
\be  \ln\frac{1+h}{1-h}\Big  |^t_{t_0}\ge 2\int_{t_0}^t f(|x_0|+t) dt,\quad t\ge t_0.\ee
We have
\begin{eqnarray}
 h(t)&&\ge 1-\frac{2}{1+\frac{1+h(t_0)}{1-h(t_0)}e^{2\int_0^t f(|x_0|+z)dz}}\nonumber\\
&& \ge 1-\frac{2}{1+e^{2\int_0^t f(|x_0|+z)dz}},\quad t\ge t_0.
\end{eqnarray}
Hence, for $T\ge c(x_0)$,
  \begin{eqnarray}\label{cmp.44}
|\gamma (T)| && =\int_{c(x_0)}^{T}h(t)dt+|\gamma (c(x_0))|\nonumber\\
&& \ge \int_{c(x_0)}^{T} \left(1-\frac{2}{1+e^{2\int_0^t f(|x_0|+z)dz}}\right)dt +|\gamma (c(x_0))|.
\end{eqnarray}

The estimate (\ref{cmp.34})  holds.

If
\be  \int_{r_c}^{+\infty} f(z)dz =+\infty, \ee
then
 \be \lim_{t\rightarrow +\infty} \left(1-\frac{2}{1+e^{2\int_0^t f(|x_0|+z)dz}}\right)=1.\ee
From (\ref{cmp.44}), \be  \lim_{t\rightarrow +\infty}\frac{|\gamma (t)|}{t}=1.\ee

The estimate (\ref{cmp.30}) holds. $\Box$

$$ $$

{\bf Proof of Theorem \ref{cmp.3_1} }\,\,\,

   With Lemma \ref{cmp.4}, we have
\be D^2r^2(X,X)\ge  2\rho_0|X|^2_g\quad \textmd{for }\ \ X\in \R^n_x, x\in \R^n.\ee
The estimate (\ref{wg.30_4})  follows from Lemma \ref{cmp.5}.$\Box$
$$ $$

{\bf Proof of Theorem \ref{cmp.3_2} }\,\,\,

   With Lemma \ref{cmp.40}, we have
\be  D^2r^2(X,X)\ge  2\rho_c|X|^2_g\quad \textmd{for }\ \ X\in \R^n_x,|x|< r_c,\ee
 \be \label{cmp.33}D^2 r(X,X)\geq
f(r)|X|_g^2\quad  \textmd{for }\ \ X\in S(r)_x,|x|\ge  r_c.\ee
The estimates (\ref{cmp.38}) and (\ref{cmp.39})  follow from Lemma \ref{cmp.31}.$\Box$
$$ $$

{\bf Proof of Theorem \ref{cmp.3_3} }\,\,\,
By a similar proof with the Theorem \ref{cmp.3_2}, the conclusions in Theorem \ref{cmp.3_3} hold true.
$$ $$

{\bf Proof of Theorem \ref{cmp.3_4} }\,\,\,
    With (\ref{pde.5}), we obtain
    \be  D^2r(X,X)=0\quad \textmd{for }\ \ X\in S(R_0)_x,|x|= R_0.\ee

    Let $\widehat{g}$ be a Riemannian metric induced by $g$ in $S(R_0)$ and $\widehat{D}$ be the associated Levi-Civita connection.

    Let $\widehat{\gamma}(t)$ be a unit-speed geodesic in  $(S(R_0),\widehat{g})$ starting at $x\in S(R_0)$, then
 \be \left\<{\widehat{\gamma}'(t)},\frac{\partial}{\partial r}\right\>_g=0,\quad \widehat{D}_{\widehat{\gamma}'(t)}{\widehat{\gamma}'(t)}=0.\ee
Therefore,
 \begin{eqnarray}  D_{\widehat{\gamma}'(t)}{\widehat{\gamma}'(t)}
&&=\widehat{D}_{\widehat{\gamma}'(t)}{\widehat{\gamma}'(t)}+\left\<D_{\widehat{\gamma}'(t)}{\widehat{\gamma}'(t)},\frac{\partial}{\partial r}\right\>_g\frac{\partial}{\partial r}
\nonumber\\
&& =\widehat{D}_{\widehat{\gamma}'(t)}{\widehat{\gamma}'(t)}-D^2r(\widehat{\gamma}'(t), \widehat{\gamma}'(t))\frac{\partial}{\partial r}
=0, \end{eqnarray}
which implies  $\widehat{\gamma}(t)$ is also the geodesic of $(\R^n,g)$.
Then
  \be\gamma(t)=\widehat{\gamma}(t)\in S(R_0),\quad \forall t\ge 0,\ee
for unit-speed geodesic $\gamma (t)$ of  $(\R^n,g)$ satisfying
 \be\gamma(0)=\widehat{\gamma}(0),\quad \gamma' (0)=\widehat{\gamma}' (0).\ee $\Box$

 \vskip .5cm
\def\theequation{5.\arabic{equation}}
\setcounter{equation}{0}
\section{Proofs for the Wave Equation }

\subsection{The propagation property }
\vskip .2cm


\quad \ \ {\bf Proof of Theorem  1.4  }\,\,\,

Let
\be  \widetilde{E}(t)=\frac12\int_{\Omega\setminus\Omega(R_0+t)}\left(u_t^2+|\nabla_gu|^2_g\right)dx_g.\ee
Then
\begin{eqnarray}
\widetilde{E}'(t) && =\frac12\frac{d}{dt}\int_{R_0+t}^{+\infty}dz\int_{x\in \Om,|x|=z} (u_t^2+|\nabla_g u|^2_g)dx_g\nonumber\\
&&=-\frac12\int_{|x|=R_0+t}\left(u_t^2+|\nabla_gu|^2_g\right)d\Ga_g +\int_{\Omega\setminus\Omega(R_0+t)}\div_g u_t \nabla_g udx_g
 \nonumber\\
&&=-\frac12\int_{|x|=R_0+t}\left(u_t^2+|\nabla_gu|^2_g\right)d\Ga_g -\int_{|x|=R_0+t} u_tu_r d\Ga_g
\nonumber\\
&&\le-\frac12\int_{|x|=R_0+t}\left(u_t^2+|\nabla_gu|^2_g\right)d\Ga_g  +\frac12\int_{|x|=R_0+t} \left(u_t^2+u_r^2\right)d\Ga_g\nonumber\\
&&\leq 0.
\end{eqnarray}

Noting that $\widetilde{E}(0)=0$, we have $u(x,t)=0,\quad |x|\ge R_0+t$.$\Box$

 \vskip .5cm
\subsection{Uniform decay of local energy for radial solutions }
\vskip .2cm

\quad \ \ In this chapter, we shall study the differences of the decay rate  between the constant coefficient wave equation and  the  wave equation on Riemannian manifold for radial solutions.

    Let $\Ga=\{x|x\in\R^n,|x|=r_0\}$ and let $u_0(x), u_1(x)$ be of compact support.

It is well-known  that the local energy for the constant coefficient wave equation
\begin{equation}
\label{wg.c11} \cases{u_{tt}- \Delta  u=0\qquad~~~~~~~~~ (x,t)\in \Om\times
(0,+\infty),\cr
 u|_{\Ga}=0\qquad~~~~~~~ t\in(0,+\infty),
\cr u(x,0)=u_0(x),\quad u_t(x,0)=u_1(x)\qquad x\in \Omega}
\end{equation}
has a uniform decay rate as follows: For  even dimensional space, the uniform decay rate of the local energy is polynomial and  for  odd dimensional space, the uniform decay rate of the local energy  is exponential. See \cite{w6},\cite{w7}.

Let $u_0(x)=u_0(r), u_1(x)=u_1(r)$
 and let $u(r,t)$ solve the following system
\begin{equation}
\label{wg.11} \cases{u_{tt}-u_{rr}-(\Delta_ gr) u_{r}=0\qquad (r,t)\in (r_0,+\infty)\times
(0,+\infty),\cr
 u(r_0)=0,
\cr u(0)=u_0(r),\quad u_t(0)=u_1(r)\qquad x\in \Omega.}
\end{equation}

Note that
\begin{eqnarray} &&\label{wg.11_1}\Delta_g  u =\div_g\left(u_{r}\frac{\partial}{\partial r}+\nabla_{\Ga_g}u\right)=u_{rr}+\Delta_g ru_{r} \quad for \ \ x \in \Om.\end{eqnarray}
Then $u(x,t)=u(r(x),t)$ solves the  system (\ref{wg.1}). The energy and the local energy for the system (\ref{wg.1})  can be rewritten as

 \be E(t,a)=\frac12\int_{\Omega(a)}\left(u_t^2+u_r^2\right)dx_g,\ee
\begin{equation}
 E(t)=\frac12\int_{\Omega}\left(u_t^2+u_r^2\right)dx_g.
\end{equation}

Using the conclusion of  scattering theory of the constant coefficient wave equation (\cite{w4}), for any positive integer $k$ , we have:
\begin{itemize}
\item
If $\Delta_g r=\frac{2k}{r}$ in (\ref{wg.11}), the   decay rate of the local energy for the system (\ref{wg.1}) is exponential, whether the dimension $n$ is even or odd.\item
If $\Delta_g r=\frac{2k-1}{r}$ in (\ref{wg.11}), the   decay rate of the local energy for the system (\ref{wg.1}) is polynomial, whether the dimension $n$ is even or odd.\end{itemize}

\begin{rem}
Let $g$ satisfy
\be \det(G(x)) =C_1(\theta)r^{4k-2(n-1)}, \quad |x|\ge r_0,\ee
where $C_1(\theta)$ is a positive function.
With (\ref{cmp.61}), we have
\be \Delta_ gr=\frac{n-1}{r}+\frac{\partial\ln \sqrt{\det(G(x))}}{\partial r}=\frac{2k}{r},\quad |x|\ge r_c.  \ee
\end{rem}
\begin{rem}
Let $g$ satisfy
\be \det(G(x)) =C_1(\theta)r^{4k-2n}, \quad |x|\ge r_0,\ee
where $C_2(\theta)$ is a positive function.
With (\ref{cmp.61}), we have
\be \Delta_ gr=\frac{n-1}{r}+\frac{\partial\ln \sqrt{\det(G(x))}}{\partial r}=\frac{2k-1}{r},\quad |x|\ge r_c.  \ee
\end{rem}


\vskip .5cm
\subsection{Uniform decay of local energy for general solutions}
\vskip .2cm

\begin{lem}\label{2wg.1} Let $u\in\widehat{H}^1_{0}(\Om)$ be of compact support. 
 Then
\be \label{1wg.16.2}\int_\Om uu_rdx_g=-\int_\Om\frac{m_2 }{2r}u^2 dx_g.  \ee
\end{lem}

{\bf Proof}.
With (\ref{cmp.61}), we have \be \Delta_ gr=\frac{n-1}{r}+\frac{\partial\ln \sqrt{\det(G(x))}}{\partial r}=\frac{m_2}{r},\quad |x|\ge r_0.\ee
Then for $|x|\ge  r_0$,
\be \div_g u^2 \frac{\partial}{\partial r}=2uu_r+u^2\Delta_g r=2uu_r+\frac{m_2}{r} u^2.\ee
Integrating over $\Om$, we have
\be\int_\Om uu_rdx_g=-\int_\Om\frac{m_2 }{2r}u^2 dx_g.\ee
The estimate (\ref{1wg.16.2}) holds.
$\Box$

\begin{lem}\label{2wg.3} Let $u_0,u_1$ be of compact support. Let u solve the system (\ref{wg.1}).
 Then
\be \label{1wg.20}\int_{\Om}r\left(u_t^2+\left |\nabla_g
u\right|_g^2\right) dx_g \leq \int_0^T\int_{\Om} (u_t^2+u_r^2) dx_g dt+\int_{\Om} r\left(u_1^2+\left |\nabla_g
u_0\right|_g^2\right) dx,\quad \forall t\ge 0.  \ee
\end{lem}

{\bf Proof}.
Let $Q=(r-t)$, letting $a\rightarrow+\infty$, it follows from (\ref{wg.14.3}) that
\begin{eqnarray}
 \int_{\Om}(r-t)\left(u_t^2+\left |\nabla_g
u\right|_g^2\right) dx_g \Big|^T_0
=&& -\int_0^T
\int_{\Om}\left(u_t^2+\left |\nabla_g
u\right|_g^2\right)  dx_gdt-2\int_0^T
\int_{\Om}u_tu_r dx_g dt\nonumber\\
\leq&&-\int_0^T
\int_{\Om}\left |\nabla_{\Ga_g}
u\right|_g^2 dx_gdt.\end{eqnarray}
Then
\begin{eqnarray}
 \int_{\Om}(r-T)\left(u_t^2+\left |\nabla_g
u\right|_g^2\right) dx_g &&\leq-\int_0^T
\int_{\Om}\left |\nabla_{\Ga_g}
u\right|_g^2 dx_gdt\nonumber\\
&&\quad +\int_{\Om} r\left(u_1^2+\left |\nabla_g
u_0\right|_g^2\right) dx.\end{eqnarray}
Simple calculation shows that
\begin{eqnarray}
 \int_{\Om}r\left(u_t^2+\left |\nabla_g
u\right|_g^2\right)  dx_g&&\leq TE(0)  -\int_0^T
\int_{\Om}\left |\nabla_{\Ga_g}
u\right|_g^2 dx_g dt +\int_{\Om} r\left(u_1^2+\left |\nabla_g
u_0\right|_g^2\right) dx
\nonumber\\
&&= \int_0^T\int_{\Om} (u_t^2+u_r^2) dx_g dt+\int_{\Om} r\left(u_1^2+\left |\nabla_g
u_0\right|_g^2\right) dx.\end{eqnarray}
The estimate (\ref{1wg.20}) holds.
$\Box$

\begin{lem}\label{2wg.5} Let$u_0,u_1$ be of compact support.  Let u solve the system (\ref{wg.1}).
 Then
\be \label{1wg.20.3} \int_{\Om}e^{r-t}\left(u_t^2+\left |\nabla_g
u\right|_g^2\right) dx_g \leq \int_{\Om}e^{r}\left(u_1^2+\left |\nabla_g
u_0\right|_g^2\right)  dx_g,\quad \forall t\ge 0.  \ee
\end{lem}

{\bf Proof}.
Let $Q=e^{r-t}$, letting $a\rightarrow+\infty$, it follows from (\ref{wg.14.3}) that
\begin{eqnarray}
 \int_{\Om}e^{r-t}\left(u_t^2+\left |\nabla_g
u\right|_g^2\right) dx_g \Big|^T_0=&&-\int_0^T
\int_{\Om}e^{r-t}\left(u_t^2+\left |\nabla_g
u\right|_g^2\right)  dx_gdt\nonumber\\&&-2\int_0^T
\int_{\Om}e^{r-t}u_tu_r dx_g dt\nonumber\\
\leq&&-\int_0^T
\int_{\Om}e^{r-t}\left(u_t^2+\left |\nabla_g
u\right|_g^2\right)  dx_gdt\nonumber\\&&+\int_0^T
\int_{\Om}e^{r-t}(u^2_t+u^2_r )dx_g dt\nonumber\\
\leq&& 0.\end{eqnarray}
Then
\be   \int_{\Om}e^{r-T}\left(u_t^2+\left |\nabla_g
u\right|_g^2\right) dx_g  \leq \int_{\Om}e^{r}\left(u_1^2+\left |\nabla_g
u_0\right|_g^2\right) dx_g .\ee
The estimate (\ref{1wg.20.3}) holds.
$\Box$

\begin{lem}\label{1wg.17} Let all the assumptions in Theorem \ref{cmp.3_6} hold. Let u solve the system (\ref{wg.1}).
 Then
\be \label{1wg.17.1}\int_0^T\int_{\Om}\left(u_t^2+u^2_r+(2m-1) |\nabla_{\Ga_g}u|^2_g\right) dx_g dt\leq \int_{\Om}r \left(u^2_t+ u^2_r\right)dx_g +C(R_0)E(0).  \ee
\end{lem}
{\bf Proof}.
Let $\H=r\frac{\partial}{\partial r}$ in (\ref{wg.14.1}). With (\ref{wg.7_1}), (\ref{cmp.41}) and (\ref{cmp.61}),  we deduce that
\begin{eqnarray} \label{wg.16_1}D\H\left(\nabla_gu,\nabla_g u\right)=&&D\H\left(u_r\frac{\partial}{\partial r} ,u_r\frac{\partial}{\partial r}\right)+D\H(\nabla_{\Ga_g}u,\nabla_{\Ga_g}u)\nonumber\\
=&&u^2_r+rD^2r\left(\nabla_{\Ga_g}u,\nabla_{\Ga_g}u\right)=u^2_r+m_1 |\nabla_{\Ga_g}u|^2_g,\quad |x|\ge r_0,\end{eqnarray}
and
\be\label{wg.16_2}\quad  \div \H= 1+r\Delta_ gr=1+r\left(\frac{n-1}{r}+\frac{\partial\ln \sqrt{\det(G(x))}}{\partial r}\right)=1+m_2, \quad |x|\ge r_0.\ee

From  (\ref{wg.14.1}),  we obtain
\begin{eqnarray}
 \label{wg.16}
 &&\int_0^T\int_{\partial\Om(a)}\frac{\pa u}{\pa\nu}\H(u) d\Ga_g dt+\frac12\int_0^T\int_{\partial\Om(a)}
\left(u_t^2-\left|\nabla_g u\right|_g^2\right)\<\H,\nu\>_g d\Ga_g dt\nonumber\\
 =&&(u_t,\H(u))\Big |^T_0+\int_0^T\int_{\Om(a)}\left(u^2_r+m_1 |\nabla_{\Ga_g}u|^2_g\right) dx_g dt\nonumber\\
&&+\int_0^T\int_{\Om(a)}\frac{\div \H}{2}\left(u_t^2-\left |\nabla_g
u\right|_g^2\right) dx_g dt\nonumber\\
=&& (u_t,\H(u))\Big |^T_0+\frac{1}{2}\int_0^T\int_{\Om(a)}\left(u_t^2+u^2_r+(2m_1-1) |\nabla_{\Ga_g}u|^2_g\right) dx_g dt\nonumber\\
&&+\int_0^T\int_{\Om(a)}\frac{m_2}{2}\left(u_t^2-\left |\nabla_g
u\right|_g^2\right) dx_g dt.
\end{eqnarray}
Let $P=\frac{m_2}{2}$, substituting (\ref{wg.14.2}) into (\ref{wg.16}), letting $a\rightarrow\infty$,  we derive that
 \begin{eqnarray}
 \label{wg.16.1}
 &&\int_0^T\int_{\Ga}\frac{\pa u}{\pa\nu}\H(u) d\Ga_g dt+\frac12\int_0^T\int_{\Ga}
\left(u_t^2-\left|\nabla_g u\right|_g^2\right)\<\H,\nu\>_g d\Ga_g dt\nonumber\\
&&+\int_0^T\int_{\Ga}Pu\frac{\pa
u}{\pa\nu}d\Ga_g dt-\frac12\int_0^T\int_{\Ga}u^2\frac{\pa
P}{\pa\nu}d\Ga_g dt\nonumber\\
=&&\left(u_t,ru_r+\frac{m_2}{2}u\right)\Big|^T_0+\frac{1}{2}\int_0^T\int_{\Om(a)}\left(u_t^2+u^2_r+(2m_1-1) |\nabla_{\Ga_g}u|^2_g\right) dx_g dt.
\end{eqnarray}

 Set
\begin{eqnarray}
 \label{wg.16.1.111}
\Pi_{\Ga}=&&\int_0^T\int_{\Ga}\frac{\pa u}{\pa\nu}\H(u) d\Ga_g dt+\frac12\int_0^T\int_{\Ga}
\left(u_t^2-\left|\nabla_g u\right|_g^2\right)\<\H,\nu\>_g d\Ga_g dt\nonumber\\
&&+\int_0^T\int_{\Ga}Pu\frac{\pa
u}{\pa\nu}d\Ga_g dt-\frac12\int_0^T\int_{\Ga}u^2\frac{\pa
P}{\pa\nu}d\Ga_g dt.
\end{eqnarray}
 Since $u\Large|_{\Ga}=0,$ we
obtain $\nabla_{\Ga_g} u\large|_{\Ga}=0$, that is,
 \be\nabla_g u=\frac{\pa
 u}{\pa\nu}\nu\quad\mbox{for}\quad
 x\in\Ga.\label{wg.16.1.2}\ee
Similarly, we have \be \H(u)=\<\H,\nabla_g u\>_g=\frac{\pa
 u}{\pa\nu}\<\H,\nu\>_g\quad\mbox{for}\quad
 x\in\Ga.\label{wg.16.1.3}\ee
Using the formulas (\ref{wg.16.1.2}) and (\ref{wg.16.1.3}) in the formula
(\ref{wg.16.1.111}) on the portion $\Ga$, we arrive at
\begin{equation}
\label{wg.16.1.4} \Pi_{\Ga}=\frac12\int_0^T
\int_{\Ga}\left(\frac{\pa
u}{\pa\nu}\right)^2\<\H\cdot\nu\>_g d\Ga_g dt=\frac12\int_0^T
\int_{\Ga}r\left(\frac{\pa
u}{\pa\nu}\right)^2\frac{\pa
 r}{\pa\nu} d\Ga_g dt\le 0.
\end{equation}

Substituting (\ref{wg.16.1.4}) into (\ref{wg.16.1}),we have
 \begin{eqnarray}
 \label{wg.16.2}
 &&\left(u_t,ru_r+\frac{m_2}{ 2} u\right)\Big|^T_0+\frac{1}{2}\int_0^T\int_{\Om}\left(u_t^2+u^2_r+(2m_1-1) |\nabla_{\Ga_g}u|^2_g\right) dx_g dt\le0.
\end{eqnarray}
With  (\ref{1wg.16.2}) and (\ref{wg.16.2}), we deduce that
\begin{eqnarray}
 \label{wg.20.5}
 &&\int_0^T\int_{\Om}\left(u_t^2+u^2_r+(2m_1-1) |\nabla_{\Ga_g}u|^2_g\right) dx_g dt\nonumber\\
 \leq&&  2\int_{\Om} r \left|u_t\left( u_r +\frac{m_2}{2r}u\right)\right| dx_g +C(R_0)E(0)\nonumber\\
 \leq&& \int_{\Om} r \left(u^2_t+ u^2_r +\frac{m_2}{r}uu_r+\frac{m_2^2}{4r^2}u^2\right)dx_g +C(R_0)E(0)\nonumber\\
 =&& \int_{\Om}  r\left(u^2_t+ u^2_r\right)dx_g +\int_{\Om} \left(m_2 uu_r+\frac{m_2^2}{4r}u^2\right)dx_g+C(R_0)E(0) \nonumber\\
=&&  \int_{\Om}  r\left(u^2_t+ u^2_r\right)dx_g -\int_{\Om} \frac{m_2^2}{4r}u^2 dx_g+C(R_0)E(0) \nonumber\\
\leq&& \int_{\Om}  r\left(u^2_t+ u^2_r\right)dx_g+C(R_0)E(0) .
\end{eqnarray}
Then, the estimate (\ref{1wg.17.1}) holds.
$\Box$
$$ $$

{\bf Proof of Theorem \ref{cmp.3_6} }\,\,\,


 Substituting (\ref{1wg.20}) into  (\ref{1wg.17.1}), we obtain
\be \label{1wg.20.1}\int_0^T\int_{\Om} |\nabla_{\Ga_g}u|^2_g dx_g dt\leq  C(R_0)E(0).  \ee
 Substituting (\ref{1wg.20.1}) into  (\ref{1wg.17.1}), we obtain
\be TE(T)\leq \int_{\Om}r(u_t^2+u_r^2)dx_g +C(R_0)E(0)\leq  \int_{\Om \backslash\Om(a) }r\left(u_t^2+\left |\nabla_g
u\right|_g^2\right)dx_g +C(a,R_0)E(0).  \ee
With (\ref{1wg.20.3}), we deduce that
\begin{eqnarray}TE(T,a)&&\leq  \int_{\Om \backslash\Om(a) }(r-T)\left(u_t^2+\left |\nabla_g
u\right|_g^2\right)dx_g +C(a,R_0)E(0)\nonumber\\
&& \leq \int_{\Om  }e^{r-T}\left(u_t^2+\left |\nabla_g
u\right|_g^2\right)dx_g +C(a,R_0)E(0)\leq C(a,R_0)E(0).  \end{eqnarray}

The estimate (\ref{wg.7.1}) holds.

\vskip .5cm
\subsection{Space-time energy estimation for general solutions}
\vskip .2cm
\begin{lem} 
Let $u\in\widehat{H}^1_{0}(\Om)$ be of compact support. Then
\be \label{1wg.16.3}\int_\Om r^{-s_2} uu_r dx_g\leq -\int_\Om\frac{m_2 r^{-2s_2} }{4}u^2 dx_g.  \ee
\end{lem}

{\bf Proof}.
With (\ref{cmp.61}), we have
\be \Delta_g r=\frac{n-1}{r}+\frac{\partial\ln \sqrt{\det(G(x))}}{\partial r}= m_2 r^{-s_2},\quad |x|\ge r_0.\ee
Then for $|x|\ge r_0$,
\begin{eqnarray} \div_g u^2 r^{-s_2}\frac{\partial}{\partial r}&&=2r^{-s_2}uu_r-s_2 r^{-s_2-1}u^2+u^2r^{-s_2}\Delta_g r\nonumber\\
&&=2r^{-s_2}uu_r+(m_2r^{-2s_2}-s_2 r^{-s_2-1}) u^2.\end{eqnarray}

Integrating over $\Om$, we have
\be2\int_\Om r^{-s_2} uu_r dx_g=-\int_\Om r^{-2s_2}(m_2 -s_2r^{s_2-1}) u^2  dx_g.\ee
Note that
\be s_2r^{s_2-1}\leq s_2r_0^{s_2-1}<\frac{(s_2+1)}{2}r_0^{s_2-1}< \frac{m_2}{2},\quad r\ge r_0.\ee
The estimate (\ref{1wg.16.3}) holds.
$\Box$

\begin{lem}\label{cmp.24} Let Assumption {\bf(B)} holds true. Let
\be P=r^{-s_2}\left(N+m_1(s_1-1)^{-1}(r_0^{1-s_1}-r^{1-s_1})\right),\quad |x|\ge r_0.\ee
 where $N$ is a positive constant. Then, for sufficiently large $N$,
\be \label{cmp.25}\Delta_g P \leq 0,\quad |x|\ge r_0. \ee
\end{lem}
{\bf Proof}.
Let
\be h(r)=\left(N+m_1(s_1-1)^{-1}(r_0^{1-s_1}-r^{1-s_1})\right),\quad r\ge r_0.\ee
Then,
\be  N\leq  h(r)\leq N+m_1(s_1-1)^{-1}r_0^{1-s_1},\quad r\ge r_0, \ee
\be  h'(r)=m_1r^{-s_1}\quad r\ge r_0.\ee

With (\ref{cmp.61}), we have
\be \Delta_g r=\frac{n-1}{r}+\frac{\partial\ln \sqrt{\det(G(x))}}{\partial r}= m_2 r^{-s_2},\quad |x|\ge r_0.\ee
We deduce that
 \begin{eqnarray}  \Delta_g P&&=\Delta_g(r^{-s_2} h(r))=\div_g\left(-s_2r^{-s_2-1}h(r)+m_1 r^{-s_1-s_2} \right)\frac{\partial}{\partial r}
\nonumber\\&&=\left(s_2(s_2+1)r^{-s_2-2}-m_2s_2 r^{-2s_2-1}\right)h(r)
 -m_1 (s_1+2s_2)r^{-s_1-s_2-1}+ m_1 m_2 r^{-s_1-2s_2} \nonumber\\
 &&\leq  s_2r^{-2s_2-1}\left((s_2+1)r^{s_2-1}-m_2\right)h(r)+ m_1 m_2 r^{-s_1-2s_2} ,\quad |x|\ge r_0. \end{eqnarray}
 Note that
\be r^{s_2-1}\leq  r_0^{s_2-1},\quad  r^{1-s_1}\leq r_0^{1-s_1}, \quad |x|\ge r_0,\ee
 \be (s_2+1)r_0^{s_2-1}-m_2<0.\ee
 Then
 \begin{eqnarray}
r^{2s_2+1}\Delta_g P&&\leq   s_2\left((s_2+1)r^{s_2-1}-m_2\right)h(r) + m_1 m_2 r^{1-s_1}
\nonumber\\
&&\leq   s_2\left((s_2+1)r_0^{s_2-1}-m_2\right)h(r) + m_1 m_2 r_0^{1-s_1}
\nonumber\\
&& \leq N s_2\left((s_2+1)r_0^{s_2-1}-m_2\right) + m_1 m_2r_0^{1-s_1},\quad |x|\ge r_0.
\end{eqnarray}
For sufficiently large $N$,
\be\Delta_g P\leq 0,\quad |x|\ge r_0.\ee
The estimate (\ref{cmp.25}) holds.
$\Box$

{\bf Proof of Theorem \ref{cmp.3_7} }\,\,\,

Let
\be h(r)=\left(N+m_1(s_1-1)^{-1}(r_0^{1-s_1}-r^{1-s_1})\right),\quad r\ge r_0.\ee
where $N\ge 1$ is a positive constant. Then,
\be  N\leq  h(r)\leq N+m_1(s_1-1)^{-1}r_0^{1-s_1},\quad r\ge r_0, \ee
\be  h'(r)=m_1r^{-s_1}\quad r\ge r_0.\ee

Let \be \H=h(r)\frac{\partial}{\partial r}\ \ in \ \ (\ref{wg.14.1}), \ee
With (\ref{cmp.9}), (\ref{cmp.41}) and (\ref{cmp.61}),  we deduce that
\begin{eqnarray} \label{cmp.13}D\H(\nabla_gu,\nabla_g u)&&=D\H(u_r\frac{\partial}{\partial r} ,u_r\frac{\partial}{\partial r})+D\H(\nabla_{\Ga_g}u,\nabla_{\Ga_g}u)\nonumber\\
&& =m_1r^{-s_1}u^2_r+h(r)D^2r(\nabla_{\Ga_g}u,\nabla_{\Ga_g}u)\nonumber\\
&&\ge m_1r^{-s_1} |\nabla_{g}u|^2_g,\quad |x|\ge r_0,\end{eqnarray}
and
\begin{eqnarray}\label{cmp.14}\div \H && = m_1r^{-s_1}+h(r)\Delta_ gr=m_1r^{-s_1}+h(r)\left(\frac{n-1}{r}+\frac{\partial\ln \sqrt{\det(G(x))}}{\partial r}\right)\nonumber\\
&&=m_1r^{-s_1}+m_2 r^{-s_2}h(r),\quad |x|\ge r_0.\end{eqnarray}

From (\ref{wg.14.1}),  we have
\begin{eqnarray}
 \label{cmp.15}
 &&\int_0^T\int_{\partial\Om(a)}\frac{\pa u}{\pa\nu}\H(u) d\Ga_g dt+\frac12\int_0^T\int_{\partial\Om(a)}
\left(u_t^2-\left|\nabla_g u\right|_g^2\right)\<\H,\nu\>_g d\Ga_g dt\nonumber\\
=&& (u_t,\H(u))\Big |^T_0+\int_0^T\int_{\Om(a)}D\H(\nabla_g
u,\nabla_g u)dx_g dt+\int_0^T\int_{\Om(a)}\frac{\div \H}{2}\left(u_t^2-\left |\nabla_g
u\right|_g^2\right) dx_g dt\nonumber\\
\ge && (u_t,\H(u))\Big |^T_0+\frac{1}{2}\int_0^T\int_{\Om(a)}m_1r^{-s_1}\left(u_t^2+\left |\nabla_g
u\right|_g^2\right) dx_g dt\nonumber\\
&&+\frac{1}{2}\int_0^T\int_{\Om(a)}m_2 r^{-s_2}h(r)\left(u_t^2-\left |\nabla_g
u\right|_g^2\right) dx_g dt.
\end{eqnarray}
Let $P=\frac{1}{2}m_2 r^{-s_2}h(r)$, substituting (\ref{wg.14.2}) into (\ref{cmp.15}), letting $a\rightarrow\infty$,  we obtain
 \begin{eqnarray}
 \label{cmp.16}
 &&\int_0^T\int_{\Ga}\frac{\pa u}{\pa\nu}\H(u) d\Ga_g dt+\frac12\int_0^T\int_{\Ga}
\left(u_t^2-\left|\nabla_g u\right|_g^2\right)\<\H,\nu\>_g d\Ga_g dt\nonumber\\
&&+\int_0^T\int_{\Ga}Pu\frac{\pa
u}{\pa\nu}d\Ga_g dt-\frac12\int_0^T\int_{\Ga}u^2\frac{\pa
P}{\pa\nu}d\Ga_g dt\nonumber\\
\geq&&\left(u_t,\H(u)+Pu\right)\Big|^T_0+\int_0^T\int_{\Om}\frac{m_1r^{-s_1}}{2}\left(u_t^2+ |\nabla_{g}u|^2_g\right) dx_g dt\nonumber\\
&&-\frac12\int_{\Om}u^2\Delta_gP dx_g dt.
\end{eqnarray}

 Set
\begin{eqnarray}
 \label{cmp.17}
\Pi_{\Ga}=&&\int_0^T\int_{\Ga}\frac{\pa u}{\pa\nu}\H(u) d\Ga_g dt+\frac12\int_0^T\int_{\Ga}
\left(u_t^2-\left|\nabla_g u\right|_g^2\right)\<\H,\nu\>_g d\Ga_g dt\nonumber\\
&&+\int_0^T\int_{\Ga}Pu\frac{\pa
u}{\pa\nu}d\Ga_g dt-\frac12\int_0^T\int_{\Ga}u^2\frac{\pa
P}{\pa\nu}d\Ga_g dt.
\end{eqnarray}
 Since $u\Large|_{\Ga}=0,$ we
obtain $\nabla_{\Ga_g} u\large|_{\Ga}=0$, that is,
 \be\nabla_g u=\frac{\pa
 u}{\pa\nu}\nu\quad\mbox{for}\quad
 x\in\Ga.\label{cmp.18}\ee
Similarly, we have \be \H(u)=\<\H,\nabla_g u\>_g=\frac{\pa
 u}{\pa\nu}\<\H,\nu\>_g\quad\mbox{for}\quad
 x\in\Ga.\label{cmp.19}\ee
Using the formulas (\ref{cmp.18}) and (\ref{cmp.19}) in the formula
(\ref{cmp.17}) on the portion $\Ga$, we obtain
\begin{equation}
\label{cmp.20} \Pi_{\Ga}=\frac12\int_0^T
\int_{\Ga}\left(\frac{\pa
u}{\pa\nu}\right)^2\<\H\cdot\nu\>_g d\Ga_g dt=\frac12\int_0^T
\int_{\Ga}h(r)\left(\frac{\pa
u}{\pa\nu}\right)^2\frac{\pa
 r}{\pa\nu} d\Ga_g dt\le 0.
\end{equation}

Substituting (\ref{cmp.20}) into (\ref{cmp.16}), for sufficiently $N$, with(\ref{cmp.25}) we deduce that
 \begin{eqnarray}
 \label{cmp.21}
 &&\left(u_t,\H(u)+Pu\right)\Big|^T_0+\int_0^T\int_{\Om}\frac{m_1r^{-s_1}}{2}\left(u_t^2+ |\nabla_{g}u|^2_g\right) dx_g dt\leq 0.
\end{eqnarray}

With  (\ref{1wg.16.3}), we obtain
\begin{eqnarray}
 \label{cmp.23}
 &&\int_0^T\int_{\Om}m_1r^{-s_1}\left(u_t^2+ |\nabla_{g}u|^2_g\right) dx_g dt\nonumber\\
\leq&&  2\int_{\Om} h(r) \left|u_t\left( u_r +\frac{m_2}{2r^{s_2}}u\right)\right| dx_g +C(R_0)E(0)\nonumber\\
\leq&&   C\int_{\Om}\left|u_t( u_r +\frac{m_2r^{-s_2}}{2}u)\right| dx_g +C(R_0)E(0)\nonumber\\
\leq&&  C \int_{\Om} \left(u^2_t+ u^2_r +m_2r^{-s_2}uu_r+\frac{m^2_2r^{-2s_2}}{4}u^2\right)dx_g +C(R_0)E(0)\nonumber\\
 =&& C\int_{\Om}  \left(u^2_t+ u^2_r\right)dx_g +m_2C\int_{\Om}\left(r^{-s_2} uu_r+\frac{m_2r^{-2s_2}}{4}u^2\right)dx_g+C(R_0)E(0) \nonumber\\
\leq&& C\int_{\Om}  \left(u^2_t+ u^2_r\right)dx_g+C(R_0)E(0) \nonumber\\
\leq &&C(R_0)E(0).
\end{eqnarray}

The estimate (\ref{cmp.12}) holds.
$\Box$
\subsection*{Acknowledgements }
This work is supported by the National Science Foundation of China, grants  no.61473126 and no.61573342 and the Key Research Program of Frontier Sciences, Chinese Academy of Sciences, no. QYZDJ-SSW-SYS011.

\end{document}